\documentclass[12pt]{article}
\usepackage{amssymb,latexsym,
amsfonts,amsbsy}
\usepackage{amsmath,amsthm}

\newtheorem{dfe}{Definition}

\newtheorem{theo}[dfe]{Theorem}

\theoremstyle{remark}
\newtheorem{Rem}{Remark}

\makeatletter

\@addtoreset{equation}{section}
\makeatother



\begin{title}
{A new Euclidean tight 6-design}
\end{title}
\author{Eiichi Bannai, \quad  Etsuko Bannai, \quad 
Junichi Shigezumi\\
}

\begin{document}
\maketitle

\begin{abstract}Euclidean $t$-designs, which are finite weighted subsets of Euclidean space,
were defined by Neumaier-Seidel (1988).
A tight $t$-design is defined as a $t$-design
whose cardinality is equal to the known natural lower bound.

In this paper, we give a new Euclidean tight 6-design in $\mathbb{R}^{22}$.
Furthermore, we also show its uniqueness up to similar transformation  fixing the origin.

This design has the structure of coherent configuration, which was defined by Higman,
and is obtained from the properties of general permutation groups.
We also show that the design is obtained by combining two orbits of McLaughlin simple group.
\end{abstract} 

\section{Introduction}
Euclidean $t$-designs were defined by  
Neumaier-Seidel \cite{N-S} as a two step generalization
of spherical designs (cf. \cite{D-G-S}). (We note that similar concepts as 
Euclidean $t$-designs have existed  
in numerical analysis as certain cubature formulas, and in statistics 
as rotatable designs (cf. \cite{B-B-H-S}).

First we give some notation.
Let $(X, w)$ be a finite weighted subset in Euclidean space $\mathbb{R}^n$,
where $X$ is a finite subset and $w$ is a positive real valued weight function on $X$.
We say $X$ is supported by $p$ concentric spheres.
That is, there are distinct nonnegative integers $r_1$, \ldots, $r_p$, 
sphere $S_i$ of radius $r_i$ centered at the origin
and subset $X_i = X \cap S_i$ for each $1 \leq i \leq p$,
such that $X = X_1 \cup \cdots \cup X_p$.
Let $w(X_i) = \sum_{x \in X_i} w(x)$ and $|S_i| = \int_{S_i} d \sigma_i (x)$.
We denote by $\mathcal{P}(\mathbb{R}^n)$ the vector space of polynomials
in $n$ variables $x_1, \ldots, x_n$ over the fields $\mathbb{R}$ of real numbers.
Let $\text{Hom}_l(\mathbb{R}^n)$ be the subspace of $\mathcal{P}(\mathbb{R}^n)$
which consists of homogeneous polynomials of degree $l$,
and let $\mathcal{P}_l (\mathbb{R}^n) = \oplus_{i=0}^l \text{Hom}_l(\mathbb{R}^n)$.

\begin{dfe}[Euclidean $t$-design]{\upshape (see \cite{N-S})}
Let $t$ be a positive integer.
A weighted finite set $(X, w)$ in $\mathbb{R}^n$ is a Euclidean $t$-design,
if the following equation
\begin{equation}
\sum_{i=1}^p \frac{w(X_i)}{|S_i|} \int_{\boldsymbol x \in S_i} f(\boldsymbol x) \, d \sigma_i (\boldsymbol x)
 = \sum_{\boldsymbol x \in X} w (\boldsymbol u) f (\boldsymbol u)
\end{equation}
is satisfied for any polynomial $f \in \mathcal{P}_t (\mathbb{R}^n)$.
\end{dfe}

There is known a natural lower bound
for the cardinalities of Euclidean $t$-designs
 (see \cite{M,D-S,N-S,B-1,B-B-5}), and tight $t$-design 
is defined as a $t$-design whose cardinality is equal to this lower bound, 
(see \cite{D-S,N-S,B-B-3,B-B-5,B-1,B-B-H-S}).
This lower bound for $t$-design for even $t$ is straightforward, but the lower 
bound for odd $t$ is somewhat delicate. (See \cite{B-B-H-S}, or the papers 
referred there, in particular those by M\"oller \cite{M-0,M}, etc., for 
more details.)

Here, we give only for the case where $t$ is even:

\begin{theo} {\upshape (see \cite{M,D-S})}
Let $(X, w)$ be a Euclidean $2 e$-design supported by $p$ concentric spheres $S$ in $\mathbb{R}^n$.
Then
\begin{equation*}
|X| \geq \dim (\mathcal{P}_e (S))
\end{equation*}
holds, where $\mathcal{P}_e (S) = \{ f|_s : f \in \mathcal{P}_e (\mathbb{R}^n) \}$
\end{theo}

\begin{dfe} {\upshape (see \cite{D-S,B-B-3})}
\begin{itemize}
\item[$(1)$] Definitions and notations are the same as above.
If the equality $(\ref{lowerbound})$ hold,
then $(X, w)$ is called a tight $2 e$-design on $p$ concentric spheres.
\item[$(2)$] Moreover, if $\dim(\mathcal{P}_e (S)) = \dim(\mathcal{P}_e (\mathbb{R}^n)) \left( = \binom{n+e}{e} \right)$ holds,
then $(X, w)$ is called a Euclidean tight $2 e$-design of $\mathbb{R}^n$.
\end{itemize}
\end{dfe}

Many tight Euclidean $t$-designs 
have been constructed (see \cite{Baj-1,Baj-2, B-B-1,
B-B-3,B-B-H-S,B-1,B-2,H-S}). However, the studies have been so far 
mostly limited to either (i) for $n=2$ or (ii) for $t\leq 5$ or $t=7$ (for $n\geq 3$).\\

In the previous paper \cite{B-B-7},  
we observed that some of the known tight 
Euclidean $t$-designs have the structure of
coherent configuration. (In particular, tight $t$-designs 
on two concentric spheres always have this property.)  
Here note that coherent configuration is a
concept defined by Higman \cite{H-1, H-2} as a generalization of
association schemes. In \cite{B-B-7} we tried to classify certain Euclidean 
$t$-designs which have the structure of coherent configuration. 
Furthermore, trying to generalize the work of \cite{B-B-7},  
we started to study Euclidean tight $6$-designs on two concentric 
spheres with one layer (fiber) being a spherical tight $4$-design. 
(Such a tight $6$-design on two concentric spheres is automatically 
a tight $6$-design of $\mathbb R^n$ as well.)  This classification problem is not yet completed, but we were able to show that 
there is only one feasible parameter set remains if $n$ is small, 
say $3\leq n\leq 438.$ (This is the parameter set described in Section 2 
of this paper.) 
At first we were a bit surprised and excited at finding this new feasible parameter set. 
Then we noticed that there 
actually exists such a Euclidean tight $6$-design, 
by combining two orbits of the McLaughlin simple group acting 
as orthogonal transformations on the  
Euclidean space  $\mathbb R^{22}.$ 
Although these two permutation representations themselves 
are well known, 
it seems new to observe that they 
actually lead to a Euclidean tight $6$-design. (Compare this with the 
well known fact that there exists no spherical tight 
$6$-design for $n\geq 3.$) 
The main purpose of this short note is 
to describe this new design. Namely, we obtain: \\

\noindent
\begin{theo} There exists a Euclidean tight $6$-design of $\mathbb R^{22}$ supported by two concentric spheres of cardinality ${22+3\choose 3}$ with the ratio $r_2/r_1=\sqrt{11}$ of the radii $r_1$ and $r_2$, and the ratio $w_2/w_1=1/729$ of the two weights
$w_1$ and $w_2.$ 
\end{theo}

\begin{Rem}
We describe all the parameters of the associated coherent 
configuration in Section 2. Also, we remark that this Euclidean tight 
6-design is unique in $\mathbb R^{22}$ supported by two concentric spheres
with $|X_1| = 275$,
up to similar transformation  fixing the origin.
\end{Rem}

\section{Parameters of the Euclidean tight $6$-design 
in $\mathbb R^{22}$} 

Let $(X,w)$ be a tight Euclidean $6$-design of 
$\mathbb R^{22}$ supported by 2 concentric spheres of positive radii $r_1$ and $r_2.$ 
Let $X=X_1\cup X_2$, and $X_i\subset S^{21}\  (i=1,2)$. Then Lemma 1.10 in \cite{B-B-3} implies that the weight function $w$
is constant on each $X_i$. Let $w\equiv w_i$ on
$X_i$ for $i=1,2$. Also Theorem 1.5 in \cite{B-B-7}
implies that $X$ has the structure of a
coherent configuration with 2 fibers.
Also Theorem 1.8 in \cite{B-1} implies that
$X_1$ and $X_2$ are spherical $4$-designs.
In the following we assume that $X_1$ is a
tight spherical $4$-design, i.e., $|X_1|=275$. 

We define
$A(X_i,X_j)=\{\frac{\boldsymbol x\cdot\boldsymbol y}{
r_ir_j}\mid  \boldsymbol x\in X_i,\ \boldsymbol y\in X_j\
\boldsymbol x \neq\boldsymbol y\}$
for any $i,j=1,2$. 
Then (Proof of) Lemma 1.10 in \cite{B-B-3} implies that
 $|A(X_i,X_j)|\leq 3$.
 Since we must have $|X_2|=|X|-|X_1|
 = {22+3\choose 3}-275=2025$, $X_2$ must be a
 $3$-distance set, i.e., $|A(X_2,X_2)|=3$ holds.
 We can also prove that $|A(X_1,X_2)|=3$ holds. 
Also, we can prove that
\begin{eqnarray}
&&A(X_1,X_1)=\{\alpha_1=\frac{1}{6},\ \alpha_2=-\frac{1}{4}\},\\
&&A(X_2,X_2)=\{\beta_1=\frac{7}{22},
\beta_2=-\frac{1}{44},\beta_3=-\frac{4}{11}\},\\
&&A(X_1,X_2)=\{\gamma_1=\frac{1}{\sqrt{11}},\
\gamma_2=-\frac{1}{4\sqrt{11}},\ \gamma_3= -\frac{3}{2\sqrt{11}}\},
\end{eqnarray}
$\frac{w_2}{w_1}=\frac{1}{729}$, and $\frac{r_2}{r_1}=\sqrt{11}$ hold.

Define 
$\alpha_0=0,\ \beta_0=0$.
Then the structure parameters of the coherent algebra (intersection numbers of the coherent configuration)  
 $p_{\alpha_i,\alpha_j}^{\alpha_k}$ 
 $(0\leq i,j,k\leq 2)$,
 $p_{\gamma_i,\gamma_j}^{\alpha_k}$ 
 $(0\leq k\leq 2,\ 1\leq i,j\leq 3)$,
$p_{\beta_i,\beta_j}^{\beta_k}$ $(0\leq i,j,k\leq 3)$
$p_{\gamma_i,\gamma_j}^{\alpha_k}$ 
 $(0\leq k\leq 3,\ 1\leq i,j\leq 3)$,
 $p_{\alpha_i,\gamma_j}^{\gamma_k}=p_{\gamma_j,\alpha_i}^{\gamma_k}$
$(0\leq i\leq 2,\  1\leq j,k\leq 3)$,
 and
$p_{\beta_i,\gamma_j}^{\gamma_k}=p_{\gamma_j,\beta_i}^{\gamma_k}$,
$0\leq i\leq 3,\  1\leq j,k\leq 3$ are determined uniquely to the values listed below.

Until this stage we just used the values of inner product $\gamma_1,
\gamma_2,\gamma_3$ to determine the structure of the coherent configuration.
In the following we introduce the following description. 
$R_{\alpha_i}=\{(\boldsymbol x,\boldsymbol y)\in X_1\times X_1
\mid \frac{\boldsymbol x\cdot\boldsymbol y}{r_1^2}=\alpha_i\}$ for $i=0,1,2$,
$R_{\beta_i}=\{(\boldsymbol x,\boldsymbol y)\in X_2\times X_2
\mid \frac{\boldsymbol x\cdot\boldsymbol y}{r_2^2}=\beta_i\}$ for $i=0,1,2,3$,
$R_{\gamma_i^{(+)}}=\{(\boldsymbol x,\boldsymbol y)\in X_1\times X_2
\mid \frac{\boldsymbol x\cdot\boldsymbol y}{r_1 r_2}=\gamma_i\}$ for $i=1,2,3$,
and
$R_{\gamma_i^{(-)}}=\{(\boldsymbol x,\boldsymbol y)\in X_2\times X_1
\mid \frac{\boldsymbol x\cdot\boldsymbol y}{r_1 r_2}=\gamma_i\}$ for $i=1,2,3$.
Thus $X\times X$ is partitioned into $13$ subsets.
For $(\boldsymbol x,\boldsymbol y)\in R_c$ we denote 
$p_{a,b}^c=|\{\boldsymbol z\in X\mid (\boldsymbol x,\boldsymbol z)\in
R_a,\ (\boldsymbol z,\boldsymbol y)\in R_b \}|$.
We define
$13\times 13$ matrices $B_a$, intersection matrix,
whose rows and columns are indexed
by the set $\{\alpha_0,\alpha_1,\alpha_2,\beta_0,
\beta_1,\beta_2,\beta_3,\gamma_1^{(+)},\gamma_2^{(+)},\gamma_3^{(+)},
\gamma_1^{(-)},\gamma_2^{(-)},\gamma_3^{(-)}\}$ with this ordering, where 
$a\in\{\alpha_0,\alpha_1,\alpha_2,\beta_0,
\beta_1,\beta_2,\beta_3,\gamma_1^{(+)},\gamma_2^{(+)},\gamma_3^{(+)},
\gamma_1^{(-)},\gamma_2^{(-)},\gamma_3^{(-)}\}$. 
The $(b,c)$ entry of $B_a$
is defined by $B_a(b,c)=p_{a,b}^c$ for any 
$a,b,c\in \{\alpha_0,\alpha_1,\alpha_2,\beta_0,
\beta_1,\beta_2,\beta_3,\gamma_1^{(+)},\gamma_2^{(+)},\gamma_3^{(+)},
\gamma_1^{(-)},\gamma_2^{(-)}$, $\gamma_3^{(-)}\}$.

{\tiny
$$
B_{\alpha_0}=\left[
\begin{array}{ccc|cccc|ccc|ccc}
1&0&0&&&&&&&&&&\\
0&1&0&&&0&&&0&&&0&\\
0&0&1&&&&&&&&&&\\
\hline
&&&&&&&&&&&\\
&0&&&&0&&&0&&&0\\
&&&&&&&&&&&\\
\hline
&&&&&&&1&0&0&&\\
&0&&&&0&&0&1&0&&0\\
&&&&&&&0&0&1&&\\
\hline
&&&&&&&&&&&\\
&0&&&&0&&&0&&&0
\end{array}\right],\quad
B_{\alpha_1}=\left[
\begin{array}{ccc|cccc|ccc|cc}
0&1&0&&&&&&&&&\\
162&105&81&&0&&&&0&&0\\
0&56&81&&&&&&&\\
\hline
&&&&&&&&&&&\\
&0&&&0&&&&0&&0\\
&&&&&&&&&\\
&&&&&&&&&\\
\hline
&&&&&&&60&42&21&&\\
&0&&&0&&&96&105&120&0\\
&&&&&&&6&15&21&&\\
\hline
&&&&&&&&&&&\\
&0&&&&0&&&0&&0&
\end{array}\right],
$$

$$B_{\alpha_2}=\left[
\begin{array}{ccc|cccc|ccc|cc}
0&0&1&&&&&&&&\\
0&56&81&&0&&&&0&&0\\
112&56&30&&&&&&&&\\
\hline
&&&&&&&&&&\\
&0&&&0&&&&0&&0\\
&&&&&&&&&&\\
&&&&&&&&&&\\
\hline
&&&&&&&16&35&56&\\
&0&&&0&&&80&70&56&0\\
&&&&&&&16&7&0&\\
\hline
&&&&&&&&&&&\\
&0&&&&0&&&0&&0&
\end{array}\right],$$
$$B_{\beta_0}=\left[
\begin{array}{c|cccc|c|ccc}
&&&&&&&&\\
0&&0&&&0&&0\\
&&&&&&&&\\
\hline
&1&0&0&0&&&\\
0&0&1&0&0&0&&0\\
&0&0&1&0&&&\\
&0&0&0&1&&&\\
\hline
0&&0&&&0&&0\\
\hline
&&&&&&1&0&0\\
0&&0&&&0&0&1&0\\
&&&&&&0&0&1
\end{array}\right],\quad
B_{\beta_1}=\left[
\begin{array}{c|cccc|c|ccc}
0&&0&&&0&&0&\\
\hline
&0&1&0&0&&&\\
0&462&185&96&28&0&&0&\\
&0&256&291&280&&&\\
&0&20&75&154&&&\\
\hline
0&&0&&&0&&0\\
\hline
&&&&&&216&105&21\\
0&&0&&&0&240&315&336\\
&&&&&&6&42&105
\end{array}\right],
$$
$$B_{\beta_2}=\left[
\begin{array}{c|cccc|c|ccc}
0&&0&&&0&&0&\\
\hline
&0&0&1&0&&&\\
&0&256&291&280&&&\\
0&1232&776&730&784&0&&0&\\
&0&200&210&168&&&\\
\hline
0&&0&&&0&&0\\
\hline
&&&&&&320&357&336\\
0&&0&&&0&816&770&840\\
&&&&&&96&105&56
\end{array}\right],\quad
B_{\beta_3}=\left[
\begin{array}{c|cccc|c|ccc}
0&&0&&&0&&0&\\
\hline
&0&0&0&1&&&\\
&0&20&75&154&&&\\
0&0&200&210&168&0&&0&\\
&330&110&45&7&&&\\
\hline
0&&0&&&0&&0&\\
\hline
&&&&&&30&105&210\\
0&&0&&&0&240&210&120\\
&&&&&&60&15&0
\end{array}\right]$$

$$B_{\gamma_1^{(+)}}=\left[
\begin{array}{ccc|c|ccc|c}
&0&&0&&0&&0\\
\hline
&&&&1&0&0&\\
&&&&216&105&21&\\
&0&&0&320&357&336&0\\
&&&&30&105&210&\\
\hline
&0&&0&&0&&0\\
\hline
567&210&81&&&&&\\
0&336&405  &0&0&0&0&0\\
0&21&81    &&&&&
\end{array}\right],\quad
B_{\gamma_2^{(+)}}=\left[
\begin{array}{ccc|c|ccc|c}
&0&&0&&0&&0\\
\hline
&&&&0&1      &0&\\
&&&&240  &315&336&\\
&0&&0&816  &770&840&0\\
&&&&240  &210&120&\\
\hline
&0&&0&&0&&0\\
\hline
0&336&405&     && &&\\
1296&840&810&0&& &0&\\
0&120&81&  &&  &&
\end{array}\right],$$
$$
B_{\gamma_3^{(+)}}=\left[
\begin{array}{ccc|c|ccc|c}
&0&&0&&0&&0\\
\hline
&&&&0  &0&1&\\
&&&&6  &42&105&\\
&0&&0&96&105&56&0\\
&&&&60&15&0&\\
\hline
&0&&0&&0&&0\\
\hline
0    &21  &81&  &&&&\\
0    &120&81&0  &&0&&0\\
162&21  &0  &  &&&&
\end{array}\right],\quad
B_{\gamma_1^{(-)}}=\left[
\begin{array}{c|cccc|c|ccc}
&&&&&&1&0&0\\
0&&0&&&0&60&42&21\\
&&&&&&16&35&56\\
\hline
0&&0&&&0&&0\\
\hline
&77&36&20&7&&&\\
0&0&40&51&56&0&&0\\
&0&1&6&14&&&\\
\hline
0&&0&&&0&&0
\end{array}\right],$$
$$B_{\gamma_2^{(-)}}=\left[
\begin{array}{c|cccc|c|ccc}
&&&&&&0 &1&0\\
0&&&0&&0&96&105&120\\
&&&&&&80&70&56\\
\hline
0&&&0&&0&  &0&\\
\hline
&0    &40   &51  &56& &&\\
0&176&120 &110&112&0 &&0\\
&0    &16    &15 &8&  &&\\
\hline
0&&0&&&0&&0
\end{array}\right],\quad
B_{\gamma_3^{(-)}}=\left[
\begin{array}{c|cccc|c|ccc}
&&&&&&0  &0  &1\\
0&&&0&&0&6  &15&21\\
&&&&&&16&7  &0\\
\hline
0&&&0&&0&&0&\\
\hline
&0  &1  &6  &14&&&\\
0&0  &16&15&8&0&&0\\
& 22&5  &1  &0&&&\\
\hline
0&&0&&&0&&0

\end{array}\right].$$
}

\section{Proof of Theorem 1}

As is well know and described in ATLAS (page 100) and 
Wilson \cite{W}, there are two permutation representations of 
$G=McL$ (McLaughlin simple group) of degree 275 and 2025 
in which the one point stabilizers are $H_1=U_4(3)$ and $H_2=M_{22}$, 
respectively. The permutation character $\chi_1$ of the
first one is
decomposed into irreducible characters as $\chi_1=1a+22a+252a$ 
and the second $\chi_2$ as $\chi_2=1a+22a+252a+1750a.$ This implies that 
by the method described below, 
we have a coherent configurations of type 
$[3,3;4]$ in the sense of D. G. Higman \cite{H-3}, by considering the 
decomposition of the permutation characters $\chi_1$ and $\chi_2.$ 
It is easy to see, that if we take the two nonzero points $\boldsymbol x_1$ and $\boldsymbol x_2$ 
in $\mathbb R^{22}$ fixed by  $H_1$ and $H_2$ respectively, then $G$ acts on the union of the two orbits $\boldsymbol x_1^G\cup 
\boldsymbol x_2^G$. Then it is easy to show that this intransitive action form a coherent 
configuration with the same parameters as given in Section 2. 
(Here the distances of $\boldsymbol x_i$ from the origin do not affect the 
structure of coherent configurations. For example, all the 
parameters $p_{i,j}^k$ as well as 
$\gamma_1, \gamma_2, \gamma_3$ do not depend on them. 
It is easy to perform these calculations, if we use either MAGMA or GAP. (It was actually performed.) 
\par
Here we note that this is proved more theoretically, by using the known 
facts. 
In ATLAS \cite{C-S-a}, Conway-Sloane \cite{C-S}, Wilson \cite{W}, all the 275 points (corresponding 
to $G/H_1$) are explicitly described, and it is shown that the action of 
$H_2=M_{22}$ on the 275 points are very visible and divided into three orbits of lengths 22, 77 and 176. 
(See Wilson page 400.) Then, using the fact that these 275 
points form a tight spherical $4$-design, using the fundamental 
equation (cf. Venkov \cite{V}, or \cite{B-B-1}), we can easily determine 
$\gamma_1, \gamma_2, \gamma_3$, and they are identical with the 
parameters given in Section 2. This completes a proof of Theorem 1. \\

The uniqueness of the Euclidean tight $6$-design in $\mathbb R^{22}$ is proved as follows. 
Since $X_1$ forms a tight $4$-design in $S^{21}\subset \mathbb R^{22},$ and since 
the uniqueness of tight $4$-design in $S^{21}$ is known, we can fix the 275 points on 
the unit sphere. We want to determine the points on the sphere $S^{21}(r_2),$ with 
the angles to any of the 275 points in $X_1$ are one of $\{\gamma_1,\gamma_2,\gamma_3\}.$ (Actually we can fix 22 points of $X_1$, which 
are linearly independent, and then we can determine the points which have all the 
angles one of $\{\gamma_1,\gamma_2,\gamma_3\}.$
The calculation shows that there are only 4050 of them. They are divided into 
two subsets each of 2025 points correspond to the two inequivalent transitive 
permutation representations of $McL/M_{22}$ of the degree 2025, which are 
interchanged each other by an outer automorphism of $McL.$ This implies the 
uniqueness of the Euclidean tight $6$-design in $\mathbb R^{22}.$\\

\begin{Rem}
The class 3 association schemes $X_2=McL/M_{22}$ is Q-polynomial 
(but not P-polynomial), and so it is in the list of Bill Martin
(see the home page:\\ {\tt http://users.wpi.edu/\~{}martin/}). 
According to Bill Martin, he obtained this information originally from an article 
(by  A. Munemasa on "spherical designs" in the
Handbook of Combinatorial Designs, 2nd ed. Chapman and Hall/CRC, pp. 617--622.) It is interesting 
to note that $X_2$ is a spherical $4$-design.  Also, it is interesting to note that 
the characterization (uniqueness) of the association scheme $X_2$ itself by parameters 
seems not yet known at the time of this writing.
\end{Rem}\quad

\begin{Rem}\label{rem-mk}
In this section, we show that the design is obtained by combining two orbits of McLaughlin simple group,
that is $McL / U_4(3)$ and $McL / M_{22}$.
In addition, by the advise of Professor Masaaki Kitazume, we noticed that there is a {\it relationship}
between these two orbits and one orbit of the Conway group ${Co}_2$ that is ${Co}_2 / U_6(2)$.
It has $4600$ points in $\mathbb{R}^{23}$, and we can classify them into $2300$ antipodal pairs.
The {\it relationship} is one to one correspondence between these pairs and $275 + 2025$ points of two orbits.
Note that Conway group ${Co}_2$ acts on these $2300$ pairs transitively.
We can refer to Remark \ref{rem-mk1} in detail.

Note that ${Co}_2 / U_6(2)$ with $4600$ points is a spherical tight $7$-design of $\mathbb{R}^{23}$,
which is uniquely determined.
\end{Rem}

\section{Some calculations}

Here, we want to explain the method to calculate vectors of $X$.
Again by Wilson \cite{W}, we obtain vectors of $X$ from the Leech lattice $\Lambda_{24}$.

We must choose the vectors $A, B \in \Lambda_{24}$ of norm $4$ such that $(A, B) = -1$.
Then, we define the following sets of vectors
\begin{align}
X_1^0 &:= \{ x \in \Lambda_{24} \mid (x, x) = 6, \, (x, A) = 3, \, (x, B) = -3 \}, \\
X_2^0 &:= \{ x \in \Lambda_{24} \mid (x, x) = 4, \, (x, A) = 2, \, (x, B) = 0 \}.
\end{align}
In conclusion, we obtain just $275$ (resp. $2025$) vectors for $X_1^0$ (resp. $X_2^0$)
from the shell of the Leech lattice of norm $6$ (resp. $4$),
which has $16773120$ (resp. $196560$) vectors.
It may not be easy to obtain $X_1^0$ since the shell of norm $6$ has too many vectors,
but we performed this calculation.

In addition, we consider the orthogonal projection $P$
from $\mathbb{R}^{24}$ to orthogonal complement
of the space spanned by vectors $A$ and $B$.
Finally, with the adjustment on radii, we obtain
\begin{equation*}
X_1 = P(X_1^0), \quad X_2 = 3 \, P(X_2^0).
\end{equation*}
Then, we obtain the design $X = X_1 \cup X_2$.\\

In the proof of the uniqueness of the design, we obtain $4050$ vectors with
the angles to any of the 275 vectors in $X_1$ are one of $\{ \gamma_1, \gamma_2, \gamma_3 \}$.

First, we want to explain a method of calculation. We denote
\begin{equation*}
\overline{X_2} := \left\{ y \in S^{21} \; \Big| \; \forall x \in \frac{1}{r_1} X_1, \; x \cdot y \in \left\{ \frac{1}{\sqrt{11}}, \; - \frac{1}{4 \sqrt{11}}, \; - \frac{3}{2 \sqrt{11}} \right\} \right\}.
\end{equation*}
The angles between the distinct vectors of $X_1$ are in $\{ \frac{1}{6}, -\frac{1}{4} \}$.
Thus, the lattice $L$ generated by all the vectors of $\frac{2 \sqrt{3}}{r_1} X_1$ is integral.
In addition, by the definition of $\overline{X_2}$, the inner products between the vectors of $\frac{2 \sqrt{3}}{r_1} X_1$
and the vectors of $\frac{4 \sqrt{11}}{2 \sqrt{3}} \overline{X_2}$ are in $\{ 4, -1, -6 \}$.
Then, $\frac{4 \sqrt{11}}{2 \sqrt{3}} \overline{X_2}$ is included in the dual lattice $L^{\sharp}$.
Furthermore, we can choose basis $\{ e_1, \ldots, e_{22} \}$ of $L$ from $\frac{2 \sqrt{3}}{r_1} X_1$.
Then, we can also obtain dual basis $\{ {e_1}', \ldots, {e_{22}}' \}$ of $L^{\sharp}$
such that $e_i \cdot {e_j}' = \delta_{i, j}$ for every $1 \leqslant i, j \leqslant 22$.
With the dual basis, we can denote the vectors of $\frac{4 \sqrt{11}}{2 \sqrt{3}} \overline{X_2}$
by $\sum_i (5 c_i - 1) {e_i}'$ for some $c_i \in \{ 0, \pm 1 \}$.
In conclusion, we can choose vectors from the lattice which is generated by the vectors
$\{ 5 e_1, \ldots, 5 e_{22}, - \sum_i {e_i}' \}$.
This lattice is much smaller than the dual lattice $L^{\sharp}$, and the norm of the every vector is equal to $\frac{44}{3}$.
Thus, in this method, we may calculate the vectors of $\overline{X_2}$ more easily.

Now, we can obtain $4050$ vectors, half of which are the vectors of $\frac{1}{r_2} X_2$.
We denote by ${X_2}' := r_2 \overline{X_2} \setminus X_2$.
Then, the calculation shows that we have ${X_2}' = P({X_2^0}')$ where
\begin{equation*}
{X_2^0}' := \{ x \in \lambda_{24} \mid (x, x) = 4, \, (x, A) = 0, \, (x, B) = -2 \}.
\end{equation*}
Moreover, we can write
\begin{align*}
X_1^0 &= \{ x \in \lambda_{24} \mid (x, x) = 6, \, (x, (-B)) = 3, \, (x, (-A)) = -3 \},\\
{X_2^0}' &= \{ x \in \lambda_{24} \mid (x, x) = 4, \, (x, (-B)) = 2, \, (x, (-A)) = 0 \}.
\end{align*}
Comparing with the definitions of $X_1^0$ and $X_2^0$,
it is clear that $X_1 \cup {X_2}'$ is isometric to $X_1 \cup X_2$.\\

\begin{Rem}\label{rem-mk1}
By Remark \ref{rem-mk}, we can calculate $X_1$ by another method. We define
\begin{align*}
Y_{+1}^0 &:= \{ x \in \Lambda_{24} \mid (x, x) = 4, \, (x, A) = 2, \, (x, B) = 1 \}\\
Y_{+2}^0 &:= \{ x \in \Lambda_{24} \mid (x, x) = 4, \, (x, A) = 2, \, (x, B) = 0 \} \quad (= X_2^0)\\
Y_{-2}^0 &:= \{ x \in \Lambda_{24} \mid (x, x) = 4, \, (x, A) = 2, \, (x, B) = -1 \}\\
Y_{-1}^0 &:= \{ x \in \Lambda_{24} \mid (x, x) = 4, \, (x, A) = 2, \, (x, B) = -2 \}
\end{align*}
where vectors $A$, $B$ are same as above.
Then, we obtain $X_1 = P(Y_{+1}^0)$ by the same orthogonal projection $P$.

Moreover, we consider another orthogonal projection $P_0$
from $\mathbb{R}^{24}$ to orthogonal complement
of the space spanned by the only vector $A$.
We denote by $Y_{\pm i} = P_0(Y_{\pm i}^0)$ for $i \in \{ 1, 2 \}$.
Then, we obtain $4600$ vectors in $Y_{+1} \cup Y_{+2} \cup Y_{-1} \cup Y_{-2} = {Co}_2 / U_6(2)$.
Since $Y_{+i} = - Y_{-i}$ for $i \in \{ 1, 2 \}$,
we have $2300$ antipodal pairs.

Oppositely, let $X = X_1 \cup X_2$ be the Euclidean tight $6$-design, having $|X_1| = 275$ and $|X_2| = 2025$.
And let $Z_{\pm 1}$ and $Z_{\pm 2}$ be subsets on $S^{22} \subset \mathbb{R}^{23}$ defined below:
\begin{align*}
Z_{+1} &= \{ (a_1 x_1, b_1) \mid x_1 \in X_1 \},
&Z_{-1} &= \{ - (a_1 x_1, b_1) \mid x_1 \in X_1 \},\\
Z_{+2} &= \{ (a_2 x_2, b_2) \mid x_2 \in X_2 \},
&Z_{-2} &= \{ - (a_2 x_2, b_2) \mid x_2 \in X_2 \},
\end{align*}
where $a_1 = \frac{2}{\sqrt{5}}$, $b_1 = \frac{1}{\sqrt{5}}$, $a_2 = \frac{2}{3 \sqrt{5}}$, and $b_2 = \frac{1}{3 \sqrt{5}}$.
Then, $Z = Z_{+1} \cup Z_{+2} \cup Z_{-1} \cup Z_{-2}$ has $4600$ points, and $Z$ is a spherical tight $7$-design.
\end{Rem}\quad

\medskip\noindent
2000 Mathematics
Subject Classification.

Primary: 05E99, Secondary:
05E30, 51M04, 65D32.\\

\noindent
Key words and Phrases.

Euclidean design, spherical design, 
association scheme,
coherent configuration,
cubature formula\\

\begin{flushright}\makebox{
\begin{minipage}{2.8in}
{\sc\large Eiichi Bannai}\vspace{0.1in}\\
Faculty of Mathematics\\
Graduate School\\
Kyushu University\\
Motooka 744, Nishi-ku\\
Fukuoka, 819-0395, Japan\\
bannai@math.kyushu-u.ac.jp\vspace{1.25in}\\
\quad
\end{minipage}
\begin{minipage}{2.4in}
{\sc\large Etsuko Bannai}\vspace{0.1in}\\
Misakigaoka 2-8-21,\\
Itoshima-shi, 819-1136, Japan\\
et-ban@rc4.so-net.ne.jp\vspace{0.5in}\\
{\sc\large Junichi Shigezumi}\vspace{0.1in}\\
Faculty of Mathematics\\
Graduate School\\
Kyushu University\\
Motooka 744, Nishi-ku\\
Fukuoka, 819-0395, Japan\\
j.shigezumi@math.kyushu-u.ac.jp
\end{minipage}
}\end{flushright}

\end{document}